# Improved Estimator of Finite Population Mean Using Auxiliary Attribute in Stratified Random Sampling

Hemant K. Verma, Prayas Sharma and †Rajesh Singh

Department of Statistics, Banaras Hindu University

Varanasi-221005

†Corresponding author

**Abstract**

The present study discuss the problem of estimating the finite population mean using auxiliary attribute in stratified random sampling. In this paper taking the advantage of point bi-serial correlation between the study variable and auxiliary attribute, we have improved the estimation of population mean in stratified random sampling. The expressions for Bias and Mean square error have been derived under stratified random sampling. In addition, an empirical study has been carried out to examine the merits of the proposed estimator over the existing estimators.

## 1. Introduction

In survey sampling, it is usual to make use of the auxiliary information at the estimation stage in order to improve the precision or accuracy of an estimator of unknown population parameter of interest. Ratio, product and regression methods of estimation are good examples in this context. Many authors including Upadhyaya and Singh (1999), Kadilar and Cingi (2005), Khoshnevisan et al. (2007), Singh et al. (2008a, b, c), Singh et al. (2009) and Singh and Kumar (2011) suggested estimators using known population parameters of an auxiliary variable. But there may be many practical situations when auxiliary information is not available directly but is qualitative in nature, that is, auxiliary information is available in the form of an attribute. For example the height of a person may depend on the fact that whether the person is male or female. The efficiency of a Dog may depend on the particular breed of that Dog. In these situations by taking the advantage of point bi-serial correlation between the study variable y and the auxiliary attribute ϕ along with the prior knowledge of the population parameter of auxiliary attribute, the estimators of population parameter of interest can be constructed.

Taking into consideration the point bi-serial correlation between auxiliary attribute and study variable, several authors including Naik and Gupta (1996), Singh et al. (2007), Shabbir and Gupta (2007), Singh et al. (2010), Abd-Elfattah et al. (2010), Singh and Solanki (2013), Malik and Singh (2013 a, b), Sharma et al. (2013 a, b),Verma et al. (2013) proposed improved estimators of population mean. In this paper we propose an estimator using the auxiliary attribute in stratified random sampling.

Let, the size of the population is N and it is stratified into L strata. Let $h^{th}$ stratum containing $N_h$ units, where h=1, 2, 3…, L such that $\sum_{h=1}^{L} N_h = N$. A sample of size $n_h$ is drawn using simple random sampling without replacement (SRSWOR) from the $h^{th}$ stratum such that $\sum_{h=1}^{L} n_h = n$. Let $y_{hi}$ and $\varphi_{hi}$ denote observed values of study variable Y and auxiliary attribute respectively, on the $i^{th}$ unit of the $h^{th}$ stratum, where i=1, 2, 3,….., $N_h$ and h=1, 2, 3,…, L. We note that, if the $i^{th}$ unit of the population possesses attribute φ then $\varphi_{hi}=1$ otherwise 0. Let $A = \sum_{i=1}^{N} \varphi_i$, $A_h = \sum_{i=1}^{N_h} \varphi_{hi}$,

$a = \sum_{i=1}^{n} \varphi_i$ and $a_h = \sum_{i=1}^{n_h} \varphi_{hi}$ denote the total number of units in the population, population stratum h, sample and sample stratum h possessing attribute $\varphi$.

Let $P = \dfrac{A}{N}$, $P_h = \dfrac{A_h}{N_h}$ and $p_h = \dfrac{a_h}{n_h}$ denote the proportion of units in population, $h^{th}$ population stratum and $h^{th}$ sample stratum, respectively.

The usual ratio and regression estimator when auxiliary variable is attribute, are respectively given by

$$t_1 = \bar{y}_{st} \left( \dfrac{P}{p_{st}} \right) \qquad (1..1)$$

$$t_2 = \bar{y}_{st} + b(P - p_{st}) \qquad (1.2)$$

where $\bar{y}_{st} = \sum_{h=1}^{L} W_h \bar{y}_h$ and

$b = \sum_{h=1}^{L} W_h^2 \gamma_h s_{\varphi y h} \Big/ \sum_{h=1}^{L} W_h^2 \gamma_h s_{\varphi h}^2$ is sample regression coefficient of y on x.

Bias and MSE expression of $t_1$ and $t_2$ are respectively given by

$$\text{Bias}(t_1) = \bar{Y}[V_2 - V_3] \qquad (1.3)$$

$$\text{Bias}(t_2) = \beta P[V_4 - V_5] \qquad (1.4)$$

$$\text{MSE}(t_1) = \bar{Y}^2[V_1 + V_2 - 2V_3] \qquad (1.5)$$

$$\text{MSE}(t_2) = \bar{Y}^2 V_1 + \beta^2 P^2 V_2 - 2\beta P \bar{Y} V_3 \qquad (1.6)$$

## 2. Proposed Estimator

Motivated by Solanki et al. (2013), we propose an estimator using information on attribute in stratified sampling as

$$t_s = \left[\bar{y}_{st} + b(P - p_{st})\right]\left(\frac{\delta P + (1-\delta)p_{st}}{\delta p_{st} + (1-\delta)P}\right) \qquad (2.1)$$

where $\delta$ is constant and b is regression coefficient.

To obtain the Bias and Mean square errors expression of (2.1), we use large sample approximation

$$\bar{y}_{st} = \bar{Y}(1+e_0), \quad p_{st} = P(1+e_1), \quad s^2_{\varphi h} = S^2_{\varphi h}(1+e_{2h}), \quad s^2_{\varphi yh} = S^2_{\varphi yh}(1+e_{3h}), \quad b = \beta(1+e_3)(1+e_2)^{-1},$$

$$e_2 = \frac{\sum_{h=1}^{L} W_h^2 \gamma_h S^2_{\varphi h} e_{2h}}{\sum_{h=1}^{L} W_h^2 \gamma_h S^2_{\varphi h}}, \qquad e_3 = \frac{\sum_{h=1}^{L} W_h^2 \gamma_h S_{\varphi yh} e_{3h}}{\sum_{h=1}^{L} W_h^2 \gamma_h S_{\varphi yh}}$$

such that $E(e_0) = E(e_1) = E(e_2) = E(e_3) = 0$.

$$E(e_0^2) = \frac{\sum_{h=1}^{L} W_h^2 \gamma_h S^2_{yh}}{\bar{Y}^2} = V_1, \qquad E(e_1^2) = \frac{\sum_{h=1}^{L} W_h^2 \gamma_h S^2_{\varphi h}}{P^2} = V_2,$$

$$E(e_0 e_1) = \frac{\sum_{h=1}^{L} W_h^2 \gamma_h S_{\varphi yh}}{\bar{Y}P} = V_3, \qquad E(e_1 e_2) = \frac{\sum_{h=1}^{L} W_h^2 \gamma_h \mu_{30h} \alpha_h}{P} = V_4,$$

$$E(e_1 e_3) = \frac{\sum_{h=1}^{L} W_h^2 \gamma_h \mu_{21h} \alpha_h^*}{P} = V_5$$

where, $S^2_{yh} = \dfrac{\sum_{i=1}^{N_h}(y_{hi} - \bar{Y}_h)}{N_h - 1}$, $\quad S_{y\varphi h} = \dfrac{\sum_{i=1}^{N_h}(y_{hi} - \bar{Y}_h)(p_{hi} - P_h)}{N_h - 1}$, $\quad \gamma_h = \dfrac{1}{n_h} - \dfrac{1}{N_h}$

$$\alpha_h = \frac{N_h^2 W_h \gamma_h}{(N_h - 1)(N_h - 2)\sum_{h=1}^{L} W_h^2 \gamma_h S^2_{\varphi h}},$$

$$\alpha_h^* = \frac{N_h^2 W_h \gamma_h}{(N_h - 1)(N_h - 2)\sum_{h=1}^{L} W_h^2 \gamma_h S_{\varphi y h}}$$

$$\mu_{rsh} = \frac{1}{N_h}\sum(y_{hi} - \overline{Y}_h)^r (p_{hi} - P_h)^s, \quad \beta = \sum_{h=1}^{L} W_h^2 \gamma_h S_{\varphi y h} \Big/ \sum_{h=1}^{L} W_h^2 \gamma_h S_{\varphi h}^2$$

Expressing (2.1) in terms of e's we have

$$t_s = \left[\overline{Y}(1+e_0) - P\beta e_1(1+e_3)(1+e_2)^{-1}\right]\left[1+(1-\delta)e_1\right]\left(1+\delta e_1\right)^{-1} \tag{2.2}$$

We assume that $|e_2| < 1$ and $|\delta e_1| < 1$ so that the term $(1+e_2)^{-1}$ and $(1+\delta e_1)^{-1}$ can be expanded.

Expanding the right hand side of (2.2) up to the first order of approximation, we have

$$t_s - \overline{Y} = \overline{Y}\left[(1-2\delta)e_1 - \delta(1-2\delta)e_1^2 + e_0 + (1-2\delta)e_0 e_1\right] - \beta P\left[e_1 + (1-2\delta)e_1^2 - e_1 e_2 + e_1 e_3\right] \tag{2.3}$$

Taking expectation of both sides of (2.3), we get the Bias expression of $t_s$ as

$$\text{Bias}(t_s) = \overline{Y}\left[(1-2\delta)V_3 - \delta(1-2\delta)V_2\right] - \beta P\left[(1-2\delta)V_2 - V_4 + V_5\right] \tag{2.4}$$

Squaring both sides of (2.3) and neglecting the terms e's having power greater than two, we have

$$(t_s - \overline{Y})^2 = \left[(1-2\delta)^2 e_1^2 + e_0^2 + 2(1-2\delta)e_0 e_1\right] + \beta^2 P^2 e_1^2 - 2\beta P \overline{Y}\left[(1-2\delta)e_1^2 + e_0 e_1\right] \tag{2.5}$$

Taking expectations of both sides of (2.4), we get the MSE of $t_s$ given by

$$\text{MSE}(t_s) = \overline{Y}\left[V_1 + (1-2\delta)^2 V_2 + 2(1-2\delta)V_3\right] + \beta^2 P^2 V_2 - 2\beta P \overline{Y}\left[(1-2\delta)V_2 + V_3\right] \tag{2.6}$$

Differentiating (2.6) with respect to $\delta$ and equating to zero, we get the optimum value of $\delta$ as

$$\delta^* = \frac{\overline{Y}(V_2 + V_3) - P\beta V_2}{2\overline{Y}V_2} \tag{2.7}$$

By substituting the optimum value of $\delta$ in (2.6) we get the min MSE of $t_s$.

## 3. Empirical Study

To illustrate the efficiency of suggested estimators in the application, we consider the data concerning the number of teachers as the study variable (*y*) and for auxiliary attribute we use number of students classifying more or less than 750, in both primary and secondary schools as auxiliary variable for 923 districts at 6 regions (as 1:Marmara 2:Agean 3:Mediterranean 4:Central Anatolia 5:Black Sea 6:East and Southeast Anatolia) in Turkey in 2007 (source: The Turkish Republic Ministry of Education). The summary statistics of the data are given in Table 2. We used Neyman allocation for allocating the samples to different strata (Cochran, 1977).

**Table 3.1: Data Statistics**

| Values | Stratum (h) | | | | | |
|---|---|---|---|---|---|---|
| | 1 | 2 | 3 | 4 | 5 | 6 |
| $N_h$ | 127 | 117 | 103 | 170 | 205 | 201 |
| $n_h$ | 31 | 21 | 29 | 38 | 22 | 39 |
| $\overline{Y}_h$ | 707.74 | 413 | 573.17 | 424.66 | 267.03 | 393.84 |
| $S_{yh}$ | 883.835 | 944.922 | 1033.467 | 810.585 | 403.654 | 711.723 |
| $S_{\varphi h}$ | 0.213 | 0.159 | 0.253 | 0.316 | 0.284 | 0.218 |
| $P_h$ | 0.952 | 0.974 | 0.932 | 0.888 | 0.912 | 0.950 |
| $S_{\varphi y h}$ | 25.267 | 9.982 | 37.453 | 44.625 | 21.04 | 18.66 |

**Table 3.2: MSE and Percentage relative efficiency of estimators**

| Estimator | MSE | PRE with respect to $\bar{y}_{st}$ |
|---|---|---|
| $\bar{y}_{st}$ | 2229.27 | 100.00 |
| $t_1$ | 2189.33 | 101.83 |
| $t_2$ | 2185.53 | 102.00 |
| $t_s$ | 2185.54 | 102.00 |

## 4. Conclusion

From theoretical discussion and empirical study we conclude that the proposed estimator $t_s$ under optimum conditions performs better than usual ratio estimator when auxiliary variable is an attribute under stratified sampling and it the minimum MSE is similar to that of usual regression estimator. The relative efficiencies and MSE of various estimators are listed in Table 3.2.